\documentclass[12pt]{article}
\usepackage[english]{babel}
\usepackage{amssymb}
\usepackage{amsfonts}
\usepackage{amsmath}

\begin{document}

\begin{center}
\bigskip {\Large On Generalized Fibonacci Quaternions \ and \
Fibonacci-Narayana Quaternions }%
\begin{equation*}
\end{equation*}

Cristina FLAUT \ and $\ $Vitalii \ SHPAKIVSKYI%
\begin{equation*}
\end{equation*}

\bigskip
\end{center}

\textbf{Abstract.} {\small In this paper, we investigate some properties of\
generalized Fibonacci quaternions and Fibonacci-Narayana quaternions in a
generalized quaternion algebra.}

\begin{equation*}
\end{equation*}



\begin{center}
\begin{equation*}
\end{equation*}
\end{center}

{\Large 0. Introduction\bigskip }%
\begin{equation*}
\end{equation*}

The Fibonacci \ numbers was introduced \ by \textit{Leonardo of Pisa
(1170-1240) }in his book \textit{Liber abbaci}, book published in 1202 AD
(see [Ko; 01], p. 1, 3). These numbers was used as a model for investigate
the growth of rabbit populations (see [Dr, Gi, Gr, Wa; 03]). The Latin name
of Leonardo \ was\textit{\ Leonardus Pisanus}, also called \textit{Leonardus
filius Bonaccii}, shortly \textit{Fibonacci}. This name is attached to the
following sequence of numbers%
\begin{equation*}
0,1,1,2,3,5,8,13,21,....,
\end{equation*}%
with the $n$th term given by the formula:%
\begin{equation*}
f_{n}=f_{n-1}+f_{n-2,}\ n\geq 2,\ 
\end{equation*}%
where $f_{0}=0,f_{1}=1.$

Fibonacci numbers was known in India before Leonardo's time and used by the
Indian authorities on metrical sciences (see [Pa; 85], p. 230). These
numbers have many properties which was studied by many authors (see [Ho;
61], [Cu; 76], [Pa; 85], [Ko; 01]).

Narayana was an outstanding Indian mathematician of the XIV century. From
him came to us the manuscript "Bidzhahanity" (incomplete), written in the
middle of the XIV century. For Narayana was interesting summation of
arithmetic series and magic squares. In the middle of the XIV century he
proved a more general summation. Using the following sums

\begin{equation*}
1+2+3+\ldots +n=S_{n}^{(1)},
\end{equation*}%
\begin{equation*}
S_{1}^{(1)}+S_{2}^{(1)}+\ldots +S_{n}^{(1)}=S_{n}^{(2)},
\end{equation*}%
\begin{equation*}
S_{1}^{(2)}+S_{2}^{(2)}+\ldots +S_{n}^{(2)}=S_{n}^{(2)},\,\ldots \,,
\end{equation*}%
Narayana calculated that

\begin{equation}
{\ S}_{n}^{(m)}{\ =}\frac{n(n+1)(n+2)\ldots (n+m)}{1\cdot 2\cdot 3\cdot
\ldots \cdot (m+1)}{\ .}  \tag{*}
\end{equation}

Narayana applied its rules to the problem of a herd of cows and heifers (see
[Yu; 61], [Si; 36], [Si; 85], [Al, Jo; 96 ]).\medskip

\textbf{Narayana problem} ([Al, Jo; 96 ]). \textit{A cow annually brings
heifers. Every heifer, beginning from the fourth year of his life also
brings heifer. How many cows and calves will be after 20 years? }

Narayana's calculation is in the following:

1) a cow within 20 years brings 20 heifers of the first generation;

2) the first heifer of the first generation brings 17 heifers second
generation, the second heifer of the first generation brings 16 heifers
second generation etc. The total in the second generation will be $%
17+16+\ldots +1=S_{17}^{(1)}$ cows and calves;

3) the first heifer of the seventeen heifers of the second generation brings
14 heifers of the third generation, the second heifer of the seventeen
heifers of the second generation brings 13 heifers of third generation, etc.
The total heifers of the first generation bring $13+12+\ldots +1=S_{13}^{(1)}
$ heads. Now, all heifers of the second generation brings $%
S_{14}^{(1)}+S_{13}^{(1)}+\ldots +S_{1}^{(1)}=S_{14}^{(2)}$ heads in the
third generation.

Similarly, Narayana calculated total number in the herd after 20 years: 
\begin{equation*}
n=1+20+S_{17}^{(1)}+S_{14}^{(2)}+\ldots +S_{2}^{(6)}.
\end{equation*}%
Using formula (*), he obtained: 
\begin{equation*}
n=1+20+\frac{17\cdot 18}{1\cdot 2}+\frac{14\cdot 15\cdot 16}{1\cdot 2\cdot 3}%
+\frac{2\cdot 3\cdot 4\cdot 5\cdot 6\cdot 7\cdot 8}{1\cdot 2\cdot 3\cdot
4\cdot 5\cdot 6\cdot 7}=2745.
\end{equation*}

This problem can be solved in the same way that Fibonacci solved its problem
about rabbits (see [Ka; 04], [Ko; 01], [Si; 36], [Si; 85]).

In the beginning of the first year was 1 cow and 1 heifer which born. That
is had 2 heads. In the beginning of the second year and in the beginning of
the third year the number of heads increased by one. Therefore the number of
heads are 3 and 4, respectively. From the fourth year, the number of heads
in the herd is defined by recurrence formulae:%
\begin{equation*}
x_{4}=x_{3}+x_{1},x_{5}=x_{4}+x_{2},\ldots ,x_{n}=x_{n-1}+x_{n-3},
\end{equation*}%
since the number of cows for any year is equal with the number of cows of
the previous year plus the number of heifers which was born (= number of
heads that were three years ago) (see [Al, Jo; 96 ]).

We have the sequence 
\begin{equation*}
2,3,4,6,9,\ldots ,\,u_{n+1}=u_{n}+u_{n-2}.
\end{equation*}%
Computing, we obtain that $u_{20}=2745$ (see [Ka; 04], [Ko; 01], [Si; 36],
[Si; 85], [Al, Jo; 96 ]).

Now, we can consider the sequence 
\begin{equation*}
1,1,1,2,3,4,6,9,\ldots ,\,u_{n+1}=u_{n}+u_{n-2},
\end{equation*}%
with $n\geq 2,u_{0}=0,u_{1}=1,u_{2}=1.$ These numbers are called the \textit{%
Fibonacci-Narayana numbers} (see [Di, St; 03]).

In the same paper [Di, St; 03], authors proved some basic properties of
Fibonacci-Narayana numbers, namely:

1) $u_{1}+u_{2}+\ldots +u_{n}=u_{n+3}-1.$

2) $u_{1}+u_{4}+u_{7}+\ldots +u_{3n-2}=u_{3n-1}.$

3) $u_{2}+u_{5}+u_{8}+\ldots +u_{3n-1}=u_{3n}.$

4) $u_{3}+u_{6}+u_{9}+\ldots +u_{3n}=u_{3n+1}-1.$

5) $u_{n+m}=u_{n-1}u_{m+2}+u_{n-2}u_{m}+u_{n-3}u_{m+1}.$

6) $u_{2n}=u_{n+1}^{2}+u_{n-1}^{2}-u_{n-2}^{2}.${\Large \ }

7) If in the sequences $\{u_{n}\},\,n=7k+4,\,n=7k+6,\,n=7k,$ when $%
k=0,1,2,\ldots $, then $u_{n}$ is even.

8) If in the sequences $\{u_{n}\}$ \thinspace\ $n=8k,\,n=8k-1,\,n=8k-3,$
when $k=0,1,2,\ldots $, then $3\mid u_{n}.$ \newline
Another property of Fibonacci-Narayana numbers was proved in [Sh; 06]. For
all natural $n\geq 2,$ we have

\begin{equation*}
u_{n}=\sum\limits_{m=0}^{[n/3]}\complement _{\lbrack n/3]}^{m}u_{n-[n/3]-2m},
\end{equation*}%
where $[a]$ is an integer part of $a$ and $\complement _{n}^{k}=\frac{n!}{%
k!\left( n-k\right) !},~k!=1\cdot 2\cdot 3\cdot ...\cdot k,k\in \mathbb{N}$%
.\medskip

Let $\mathbb{H}\left( \beta _{1},\beta _{2}\right) $ be the generalized
real\ quaternion algebra, the algebra of the elements of the form $%
a=a_{1}\cdot 1+a_{2}e_{2}+a_{3}e_{3}+a_{4}e_{4},$ where $a_{i}\in \mathbb{R}%
, i\in \{1,2,3,4\}$, and the basis elements $1,e_2,e_3,e_4$ satisfy the
following multiplication table: \vspace{3mm}

\begin{center}
\begin{tabular}{ccccc}
$\cdot\:\,\,$\vline & $1$ & $e_2$ & $e_3$ & $e_4$ \\ \hline
$1$\:\,\vline & $1$ & $e_2$ & $e_3$ & $e_4$ \\ 
$e_2$\vline & $e_2$ & $-\beta_1$ & $e_4$ & $-\beta_1e_3$ \\ 
$e_3$\vline & $e_3$ & $-e_4$ & $-\beta_2$ & $\beta_2e_2$ \\ 
$e_4$\vline & $e_4$ & $\beta_1e_3$ & $-\beta_2e_2$ & $-\beta_1\beta_2$%
\end{tabular}
\end{center}

\vspace{3mm}

We denote by $\boldsymbol{t}\left( a\right) $ and $\boldsymbol{n}\left(
a\right) $ the trace and the norm of a real quaternion $a.$ The norm of a
generalized quaternion has the following expression $\boldsymbol{n}\left(
a\right) =a_{1}^{2}+\beta _{1}a_{2}^{2}+\beta_{2}a_{3}^{2}+\beta
_{1}\beta_2a_{4}^{2}.$ For $\beta _{1}=\beta _{2}=1,$ we obtain the real
division algebra $\mathbb{H}$.

\begin{equation*}
\end{equation*}

{\Large 1. Preliminaries}%
\begin{equation*}
\end{equation*}

In the present days, several mathematicians studied properties of the
Fibonacci sequence. In [Ho; 61], the \ author generalized Fibonacci numbers
and gave many properties of them: 
\begin{equation*}
h_{n}=h_{n-1}+h_{n-2},\ \ n\geq 2,
\end{equation*}%
where $h_{0}=p,h_{1}=q,$ with $\ p,q$ being arbitrary integers. In the same
paper [Ho; 61], relation (7), the following relation between Fibonacci
numbers and generalized Fibonacci numbers was obtained: 
\begin{equation}
h_{n+1}=pf_{n}+qf_{n+1}.  \tag{1.1}
\end{equation}%
The same author, in [Ho; 63], defined \ and studied Fibonacci quaternions
and generalized Fibonacci quaternions in the real division quaternion
algebra and found \ a lot of properties of them. For the generalized real\
quaternion algebra, the Fibonacci quaternions and generalized Fibonacci
quaternions are defined in the same way:%
\begin{equation*}
F_{n}=f_{n}\cdot 1+f_{n+1}e_{2}+f_{n+2}e_{3}+f_{n+3}e_{4},
\end{equation*}%
for the \ $n$th Fibonacci quaternions, and 
\begin{equation*}
H_{n}=h_{n}\cdot 1+h_{n+1}e_{2}+h_{n+2}e_{3}+h_{n+3}e_{4},
\end{equation*}%
for the \ $n$th generalized \ Fibonacci quaternions.

In the same paper, we find the norm formula \ for the \ $n$th Fibonacci
quaternions:

\begin{equation*}
\boldsymbol{n}\left( F_{n}\right) =F_{n}\overline{F}_{n}=3f_{2n+3}, 
\tag{1.2}
\end{equation*}%
where \ $\overline{F}_{n}=f_{n}\cdot
1-f_{n+1}e_{2}-f_{n+2}e_{3}-f_{n+3}e_{4} $ is the conjugate of the $F_{n}$
in the algebra $\mathbb{H}.$ After that, many authors studied Fibonacci and
generalized Fibonacci quaternions in the real division quaternion algebra
giving more and surprising new properties (for example, see [Sw; 73],
[Sa-Mu; 82] and [Ha; 12]).

M. N. S. Swamy, in [Sw; 73], formula (17), \ obtained the norm formula for \
the \ $n$th generalized \ Fibonacci quaternions: 
\begin{eqnarray*}
\boldsymbol{n}\left( H_{n}\right) &=&H_{n}\overline{H}_{n}= \\
&=&3(2pq-p^{2})f_{2n+2}+(p^{2}+q^{2})f_{2n+3},
\end{eqnarray*}%
where \ $\overline{H}_{n}=h_{n}\cdot
1-h_{n+1}e_{2}-h_{n+2}e_{3}-h_{n+3}e_{4} $ is the conjugate of the $H_{n}$
in the algebra $\mathbb{H}.~\ $

Similar to A.~F.~Horadam, we define the Fibonacci-Narayana quaternions as 
\begin{equation*}
U_{n}=u_{n}\cdot 1+u_{n+1}e_{2}+u_{n+2}e_{3}+u_{n+3}e_{4},
\end{equation*}%
where $u_{n}$ are the $n$th Fibonacci-Narayana number.

In this paper, we \ give some properties of \ generalized Fibonacci
quaternions and \ Fibonacci-Narayana quaternions.\ 

\bigskip

\begin{equation*}
\end{equation*}%
\begin{equation*}
\end{equation*}%
\begin{equation*}
\end{equation*}

{\Large 2. \ Generalized Fibonacci Quaternions\bigskip }%
\begin{equation*}
\end{equation*}

\qquad\ \ \ \ \ As in the case of Fibonacci numbers, numerous results
between Fibonacci generalized numbers can be deduced. In the following, we
will study some properties of the generalized Fibonacci quaternions in the
generalized real quaternion algebra $\mathbb{H}\left( \beta _{1},\beta
_{2}\right) $. Let $\ $ $F_{n}=f_{n}\cdot
1+f_{n+1}e_{2}+f_{n+2}e_{3}+f_{n+3}e_{4}$ be the $n$th Fibonacci \
quaternion and \ $H_{n}=h_{n}\cdot 1+h_{n+1}e_{2}+h_{n+2}e_{3}+h_{n+3}e_{4}$
be the $n$th generalized Fibonacci \ quaternion. A first question which can
arise is what algebraic structure have these elements? The answer will be
found in the below theorem, denoting first a $n$th generalized Fibonacci
number \ and a $n$th generalized Fibonacci element with $h_{n}^{p,q},$
respectively $H_{n}^{p,q}.$ In this way, we emphasis the starting integers $%
p $ and $q.\medskip $

\textbf{Theorem} \textbf{2.1.} \textit{The set} $\mathcal{H}%
_{n}=\{H_{n}^{p,q}~/~p,q\in \mathbb{Z}\}\cup \{0\}$ \textit{is a} $\mathbb{Z}%
-$\textit{module}.\medskip

\textbf{Proof.} Indeed, $aH_{n}^{p,q}+bH_{n}^{p^{\prime },q\prime
}=H_{n}^{ap+bp^{\prime },aq+bq^{\prime }}\in \mathcal{H}_{n},$ where\newline
$a,b,p,q,p^{\prime },q^{\prime }\in \mathbb{Z}.\Box \medskip $

\textbf{Theorem \ 2.2. } \textit{i)} \textit{For the Fibonacci quaternion
elements, we have}

\begin{equation}
\sum\limits_{m=1}^{n}\left( -1\right) ^{m+1}F_{m}=\left( -1\right)
^{n+1}F_{n-1}+1+e_{3}+e_{4}.  \tag{2.1}
\end{equation}

\textit{ii) For the generalized Fibonacci quaternion elements, the following
relation is true}%
\begin{equation}
\sum\limits_{m=1}^{n}\left( -1\right) ^{m+1}H_{m}^{p,q}\text{=}\left(
-1\right) ^{n+1}H_{n-1}^{p,q}-p\text{+}q\text{+}pe_{2}\text{+}qe_{3}\text{+}%
pe_{4}\text{+}qe_{4}  \tag{2.2}
\end{equation}

\textbf{Proof.}

i) From [Cu; 76], we know that 
\begin{equation}
\sum\limits_{m=1}^{n}\left( -1\right) ^{m+1}f_{m}=\left( -1\right)
^{n+1}f_{n-1}+1.  \tag{2.3}
\end{equation}%
It results:\newline
$\sum\limits_{m=1}^{n}\left( -1\right) ^{m+1}F_{m}$=\newline
=$\sum\limits_{m=1}^{n}\left( -1\right) ^{m+1}f_{m}$+$e_{2}\sum%
\limits_{m=1}^{n}\left( -1\right) ^{m+1}f_{m+1}$+\newline
+$e_{3}\sum\limits_{m=1}^{n}\left( -1\right) ^{m+1}f_{m+2}$+$%
e_{4}\sum\limits_{m=1}^{n}\left( -1\right) ^{m+1}f_{m+3}$= \newline
=$(-1)^{n+1}f_{n-1}+1-e_{2}[(-1)^{n+1}f_{n-1}+\left( -1\right)
^{n+2}f_{n+1}] $+\newline
+$e_{3}\left[ \left( -1\right) ^{n+1}f_{n-1}+1+\left( -1\right)
^{n}f_{n+1}+\left( -1\right) ^{n+1}f_{n+2}\right] -$\newline
$-e_{4}[\left( -1\right) ^{n+1}f_{n-1}-1+\left( -1\right)
^{n+2}f_{n+1}+\left( -1\right) ^{n+3}f_{n+2}+\left( -1\right) ^{n+4}f_{n+3}]$%
=\newline
=$(-1)^{n+1}f_{n-1}+1+\left( -1\right) ^{n+1}e_{2}f_{n}+e_{3}\left(
-1\right) ^{n+1}\left[ f_{n+1}+\left( -1\right) ^{n+1}\right] -$\newline
$-e_{4}\left( -1\right) ^{n+1}\left[ -f_{n+2}-\left( -1\right) ^{n+1}\right]
=$\newline
$=\left( -1\right) ^{n+1}\left(
f_{n-1}+f_{n}e_{2}+f_{n+1}e_{3}+f_{n+2}e_{4}\right) +1+e_{3}+e_{4}=$\newline
$=\left( -1\right) ^{n+1}F_{n-1}+1+e_{3}+e_{4}$.

ii) Using relations $\left( 1.1\right) $ and $\left( 2.3\right) ,$we have%
\newline
$\sum\limits_{m=1}^{n}\left( -1\right) ^{m+1}H_{m}^{p,q}=$\newline
$=\sum\limits_{m=1}^{n}\left( -1\right) ^{m+1}h_{m}^{p,q}$+$%
e_{2}\sum\limits_{m=1}^{n}\left( -1\right) ^{m+1}h_{m+1}^{p,q}$+\newline
+$e_{3}\sum\limits_{m=1}^{n}\left( -1\right) ^{m+1}h_{m+2}^{p,q}$+$%
e_{4}\sum\limits_{m=1}^{n}\left( -1\right) ^{m+1}h_{m+3}^{p,q}=$\newline
$=\sum\limits_{m=1}^{n}\left( -1\right)
^{m+1}pf_{m-1}+\sum\limits_{m=1}^{n}\left( -1\right) ^{m+1}qf_{m}+$\newline
$+e_{2}\sum\limits_{m=1}^{n}\left( -1\right)
^{m+1}pf_{m}+e_{2}\sum\limits_{m=1}^{n}\left( -1\right) ^{m+1}qf_{m+1}+$%
\newline
$+e_{3}\sum\limits_{m=1}^{n}\left( -1\right)
^{m+1}pf_{m+1}+e_{3}\sum\limits_{m=1}^{n}\left( -1\right) ^{m+1}qf_{m+2}+$%
\newline
$+e_{4}\sum\limits_{m=1}^{n}\left( -1\right)
^{m+1}pf_{m+2}+e_{4}\sum\limits_{m=1}^{n}\left( -1\right) ^{m+1}qf_{m+3}=$%
\newline
$=p\left( -1\right) ^{n+1}f_{n-2}-p+q\left( -1\right) ^{n+1}f_{n-1}+q+$%
\newline
$+e_{2}p\left( -1\right) ^{n+1}f_{n-1}+pe_{2}+e_{2}q\left[ \left( -1\right)
^{n+1}f_{n+1}-\left( -1\right) ^{n+1}f_{n-1}\right] +$\newline
$+e_{3}p\left[ \left( -1\right) ^{n+1}f_{n+1}-\left( -1\right) ^{n+1}f_{n-1}%
\right] +$\newline
$+e_{3}q\left[ \left( -1\right) ^{n+1}f_{n-1}+1+\left( -1\right)
^{n}f_{n+1}+\left( -1\right) ^{n+1}f_{n+2}\right] +$\newline
$+e_{4}p\left[ \left( -1\right) ^{n+1}f_{n-1}+1+\left( -1\right)
^{n}f_{n+1}+\left( -1\right) ^{n+1}f_{n+2}\right] -$\newline
$-e_{4}q[\left( -1\right) ^{n+1}f_{n-1}-1+\left( -1\right)
^{n+2}f_{n+1}+\left( -1\right) ^{n+3}f_{n+2}+\left( -1\right)
^{n+4}f_{n+3}]= $\newline
$=p\left( -1\right) ^{n+1}f_{n-2}-p+q\left( -1\right) ^{n+1}f_{n-1}+q+$%
\newline
$+e_{2}p\left( -1\right) ^{n+1}f_{n-1}+pe_{2}+e_{2}q\left( -1\right)
^{n+1}f_{n}+e_{3}p\left( -1\right) ^{n+1}f_{n}+$\newline
$+e_{3}q\left( -1\right) ^{n+1}[f_{n-1}+\left( -1\right)
^{n+1}-f_{n+1}+f_{n+2}]+$\newline
$+e_{4}p\left( -1\right) ^{n+1}\left[ f_{n-1}+\left( -1\right)
^{n+1}-f_{n+1}+f_{n+2}\right] -$\newline
$-e_{4}q\left( -1\right) ^{n+1}\left[ f_{n-1}-\left( -1\right)
^{n+1}-f_{n+1}+f_{n+2}-f_{n+3}\right] =$\newline
$=p\left( -1\right) ^{n+1}f_{n-2}-p+q\left( -1\right) ^{n+1}f_{n-1}+q+$%
\newline
$+e_{2}p\left( -1\right) ^{n+1}f_{n-1}+pe_{2}+e_{2}q\left( -1\right)
^{n+1}f_{n}+e_{3}p\left( -1\right) ^{n+1}f_{n}+$\newline
$+e_{3}q\left( -1\right) ^{n+1}\left[ f_{n+1}+\left( -1\right) ^{n+1}\right]
+e_{4}p\left( -1\right) ^{n+1}[\left( -1\right) ^{n+1}+f_{n+1}]-$\newline
$-e_{4}q\left( -1\right) ^{n+1}\left[ -f_{n+2}-\left( -1\right) ^{n+1}\right]
=$\newline
$=\left( -1\right) ^{n+1}H_{n-1}^{p,q}-p+q+pe_{2}+qe_{3}+pe_{4}+qe_{4}.\Box
\medskip \medskip $

From the above Theorem, we can remark that all identities valid for the
Fibonacci quaternions can be easy adapted in an approximative similar
expression for the generalized Fibonacci quaternions, if we use relation $%
\left( 1.1\right) ,$ a true relation in the both algebras $\mathbb{H}\left(
\beta _{1},\beta _{2}\right) $ and $\mathbb{H}.\medskip $

\textbf{Proposition 2.3.} \ \textit{If}~ $h_{n+1}=pf_{n}+qf_{n+1}=0,$ 
\textit{then we have: }

\begin{equation}
H_{n+1}^{2}\text{=}3\frac{q^{2}}{f_{n}^{2}}\left[
f_{2n+1}^{2}-f_{n+1}f_{n-2}f_{2n+2}\right] ,  \tag{2.4}
\end{equation}%
\textit{where} \ $H_{n+1}^{2}\in \mathbb{H}\left( \beta _{1},\beta
_{2}\right) .\medskip $

\textbf{Proof.} Since $h_{n+1}=0,$ it results that $\boldsymbol{t}\left(
H_{n+1}\right) =h_{n+1}=0,$ therefore\newline
$\boldsymbol{n}\left( H_{n+1}\right) =H_{n+1}^{2}$. From $%
h_{n}=pf_{n}+qf_{n+1}=0 $, we have \ $p=-\frac{qf_{n+1}}{f_{n}}$ and\newline
we obtain:%
\begin{equation*}
p^{2}+2pq\text{=}\frac{q^{2}f_{n+1}^{2}}{f_{n}^{2}}-2q^{2}\frac{f_{n+1}}{%
f_{n}}\text{=}-\frac{q^{2}f_{n+1}f_{n-2}}{f_{n}^{2}}
\end{equation*}%
and

\begin{equation*}
p^{2}+q^{2}\text{=}\frac{q^{2}f_{n+1}^{2}}{f_{n}^{2}}+q^{2}\text{=}q^{2}%
\frac{f_{n+1}^{2}+f_{n}^{2}}{f_{n}^{2}}=q^{2}\frac{f_{2n+1}}{f_{n}^{2}},
\end{equation*}%
since \ $f_{n+1}^{2}+f_{n}^{2}=f_{2n+1}.$

It results \newline
$\boldsymbol{n}\left( H_{n+1}\right) $=$3[(p^{2}+2pq)f_{2n+2}$+$%
(p^{2}+q^{2})f_{2n+1}]$=\newline
=$3\frac{q^{2}}{f_{n}^{2}}[-f_{n+1}f_{n-2}f_{2n+2}+f_{2n+1}^{2}].\Box
\medskip $

In the following, we will compute the norm of a Fibonacci quaternion and of
a generalized Fibonacci quaternion in the algebra $\mathbb{H}\left( \beta
_{1},\beta _{2}\right) .$

Let $F_{n}=f_{n}\cdot 1+f_{n+1}e_{2}+f_{n+2}e_{3}+f_{n+3}e_{4}$ be the $n$th
Fibonacci quaternion, then its norm is 
\begin{equation*}
\boldsymbol{n}\left( F_{n}\right) =f_{n}^{2}+\beta _{1}f_{n+1}^{2}+\beta
_{2}f_{n+2}^{2}+\beta _{1}\beta _{2}f_{n+3}^{2}.
\end{equation*}%
Using recurrence of Fibonacci numbers and relations 
\begin{equation}
f_{n}^{2}+f_{n-1}^{2}=f_{2n-1},\quad n\in \mathbb{N},  \tag{2.5}
\end{equation}%
\begin{equation}
f_{2n}=f_{n}^{2}+2f_{n}f_{n-1},\quad n\in \mathbb{N},  \tag{2.6}
\end{equation}%
from [Cu; 76], we have\newline
$\boldsymbol{n}\left( F_{n}\right) =f_{n}^{2}+\beta _{1}f_{n+1}^{2}+\beta
_{2}f_{n+2}^{2}+\beta _{1}\beta _{2}f_{n+3}^{2}=$\newline
$=f_{n}^{2}+\beta _{1}f_{n+1}^{2}+\beta _{2}\left( f_{n+2}^{2}+\beta
_{1}f_{n+3}^{2}\right) =$\newline
$=f_{2n+1}+\left( \beta _{1}-1\right) f_{n+1}^{2}+\beta _{2}\left(
f_{2n+5}+(\beta _{1}-1)f_{n+3}^{2}\right) =$\newline
$=f_{2n+1}+\beta _{2}f_{2n+5}+(\beta _{1}-1)\left( f_{n+1}^{2}+\beta
_{2}f_{n+3}^{2}\right) =$\newline
$=\left( 1+2\beta _{2}\right) f_{2n+1}+3\beta _{2}f_{2n+2}+(\beta
_{1}-1)\left( f_{n+1}^{2}+\beta _{2}f_{n+3}^{2}\right) =$\newline
$=h_{2n+2}^{1+2\beta _{2},3\beta _{2}}+(\beta _{1}-1)\left(
f_{n+1}^{2}+\beta _{2}f_{n+3}^{2}\right) =$\newline
$=h_{2n+2}^{1+2\beta _{2},3\beta _{2}}+(\beta _{1}-1)\left(
f_{2n+2}-2f_{n}f_{n+1}+\beta _{2}f_{2n+6}-2\beta _{2}f_{n+2}f_{n+3}\right) =$%
\newline
$=h_{2n+2}^{1+2\beta _{2},3\beta _{2}}+(\beta _{1}-1)[f_{2n+2}+\beta
_{2}f_{2n+6}-2\left( f_{n}f_{n+1}+\beta _{2}f_{n+2}f_{n+3}\right) ]=$\newline
$=h_{2n+2}^{1+2\beta _{2},3\beta _{2}}$+$(\beta _{1}-1)[f_{2n+2}$+$\beta
_{2}f_{2n+6}-2\left( f_{n}f_{n+1}\text{+}\beta _{2}f_{n+2}^{2}\text{+}\beta
_{2}f_{n+1}f_{n+2}\right) =$\newline
$=h_{2n+2}^{1+2\beta _{2},3\beta _{2}}$+$(\beta _{1}$-$1)[f_{2n+2}$+$\beta
_{2}f_{2n+6}$-$2\left( f_{n}f_{n+1}\text{+}\beta _{2}f_{n+2}^{2}\text{+}%
\beta f_{n+1}^{2}\text{+}\beta _{2}f_{n}f_{n+1}\right) $=\newline
$=h_{2n+2}^{1+2\beta _{2},3\beta _{2}}$+$(\beta _{1}-1)[f_{2n+2}$+$\beta
_{2}f_{2n+6}-2\left( 1+\beta _{2}\right) f_{n}f_{n+1}-2\beta _{2}f_{2n+3}]=$%
\newline
$=h_{2n+2}^{1+2\beta _{2},3\beta _{2}}$+$(\beta _{1}$-$1)[f_{2n+2}$+$\beta
_{2}f_{2n+4}$+$\beta _{2}f_{2n+3}$+$\beta _{2}f_{2n+4}-$\newline
$-2\beta _{2}f_{2n+3}$-$2\left( 1\text{+}\beta _{2}\right) f_{n}f_{n+1}]=$%
\newline
$=h_{2n+2}^{1+2\beta _{2},3\beta _{2}}$+$(\beta _{1}$-$1)[f_{2n+2}+2\beta
_{2}f_{2n+4}-\beta _{2}f_{2n+3}$-$2\left( 1\text{+}\beta _{2}\right)
f_{n}f_{n+1}]=$\newline
$=h_{2n+2}^{1+2\beta _{2},3\beta _{2}}$+$(\beta _{1}$-$1)[f_{2n+2}+2\beta
_{2}f_{2n+2}+2\beta _{2}f_{2n+3}-$\newline
$-\beta _{2}f_{2n+3}$-$2\left( 1\text{+}\beta _{2}\right) f_{n}f_{n+1}]=$%
\newline
$=h_{2n+2}^{1+2\beta _{2},3\beta _{2}}$+$(\beta _{1}$-$1)[\left( 1+2\beta
_{2}\right) f_{2n+2}+\beta _{2}f_{2n+3}$-$2\left( 1\text{+}\beta _{2}\right)
f_{n}f_{n+1}]=$\newline
$=h_{2n+2}^{1+2\beta _{2},3\beta _{2}}$+$(\beta _{1}$-$1)[h_{2n+3}^{1+2\beta
_{2},\beta _{2}}-2\left( 1\text{+}\beta _{2}\right) f_{n}f_{n+1}]=$\newline
$=h_{2n+2}^{1+2\beta _{2},3\beta _{2}}$+$(\beta _{1}$-$1)h_{2n+3}^{1+2\beta
_{2},\beta _{2}}-2(\beta _{1}$-$1)\left( 1\text{+}\beta _{2}\right)
f_{n}f_{n+1}.\medskip $

We just proved\medskip

\textbf{Theorem 2.4.} \textit{The norm of the} $n$th \textit{Fibonacci
quaternion} $F_{n}$ \textit{in a generalized quaternion algebra is}%
\begin{equation}
\boldsymbol{n}\left( F_{n}\right) \text{=}h_{2n+2}^{1+2\beta _{2},3\beta
_{2}}\text{+}(\beta _{1}\text{-}1)h_{2n+3}^{1+2\beta _{2},\beta _{2}}\text{-}%
2(\beta _{1}\text{-}1)\left( 1\text{+}\beta _{2}\right) f_{n}f_{n+1}.\medskip
\tag{2.7}
\end{equation}%
$\Box \medskip \medskip $

Using the formula $\left(2.7\right)$ and the relation $\left(1.1\right)$
when $\beta_1=\beta_2=1,$ we obtain the formula $\left(1.2\right)$.

Using the above theorem and relations $\left( 2.5\right) $ and $\left(
2.6\right) ,$ we can compute the norm of a generalized Fibonacci quaternion
in a generalized quaternion algebra. Let $H_{n}=h_{n}\cdot
1+h_{n+1}e_{2}+h_{n+2}e_{3}+h_{n+3}e_{4}$ be the $n$th generalized Fibonacci
quaternion. The norm is\medskip \newline
$\boldsymbol{n}\left( H_{n}^{p,q}\right) $=$h_{n}^{2}$+$\beta
_{1}h_{n+1}^{2} $+$\beta _{2}h_{n+2}^{2}$+$\beta _{1}\beta _{2}h_{n+3}^{2}$=%
\newline
=$\left( pf_{n-1}\text{+}qf_{n}\right) ^{2}$+$\beta _{1}\left( pf_{n}\text{+}%
qf_{n+1}\right) ^{2}$+$\beta _{2}\left( pf_{n+1}\text{+}qf_{n+2}\right)
^{2}+ $\newline
+$\beta _{1}\beta _{2}\left( pf_{n+2}\text{+}qf_{n+3}\right) ^{2}=$\newline
$=p^{2}\left( f_{n-1}^{2}+\beta _{1}f_{n}^{2}+\beta _{2}f_{n+1}^{2}+\beta
_{1}\beta _{2}f_{n+2}^{2}\right) +$\newline
$+q^{2}\left( f_{n}^{2}+\beta _{1}f_{n+1}^{2}+\beta _{2}f_{n+2}^{2}+\beta
_{1}\beta _{2}f_{n+3}^{2}\right) +$\newline
$+2pq\left( f_{n-1}f_{n}+\beta _{1}f_{n}f_{n+1}+\beta
_{2}f_{n+1}f_{n+2}+\beta _{1}\beta _{2}f_{n+3}f_{n+2}\right) =$\newline
$=p^{2}h_{2n}^{1+2\beta _{2},3\beta _{2}}$+$p^{2}(\beta _{1}$-$%
1)h_{2n+1}^{1+2\beta _{2},\beta _{2}}$-$2p^{2}(\beta _{1}$-$1)\left( 1\text{+%
}\beta _{2}\right) f_{n-1}f_{n}+$\newline
$+q^{2}h_{2n+2}^{1+2\beta _{2},3\beta _{2}}$+$q^{2}(\beta _{1}$-$%
1)h_{2n+3}^{1+2\beta _{2},\beta _{2}}$-$2q^{2}(\beta _{1}$-$1)\left( 1\text{+%
}\beta _{2}\right) f_{n}f_{n+1}+$\newline
$+2pq\left( 1\text{-}\beta _{1}\right) f_{n}f_{n-1}$+$2pq\beta _{1}f_{2n}$ +$%
2pq\beta _{2}\left( 1\text{-}\beta _{1}\right) f_{n+1}f_{n+2}$+$2pq\beta
_{1}\beta _{2}f_{2n+4}$=\newline
$=p^{2}h_{2n}^{1+2\beta _{2},3\beta _{2}}$+$p^{2}(\beta _{1}$-$%
1)h_{2n+1}^{1+2\beta _{2},\beta _{2}}+q^{2}h_{2n+2}^{1+2\beta _{2},3\beta
_{2}}$+$q^{2}(\beta _{1}$-$1)h_{2n+3}^{1+2\beta _{2},\beta _{2}}-$\newline
$-2p\left( \beta _{1}-1\right) \left( p\beta _{2}+p+q\right)
f_{n-1}f_{n}-2q^{2}(\beta _{1}$-$1)\left( 1\text{+}\beta _{2}\right)
f_{n}f_{n+1}+$\newline
$+h_{2n+1}^{2pq\beta _{1},2pq\beta _{1}\beta _{2}}+2pq\beta _{1}\beta
_{2}(f_{2n}+f_{2n+3})+2pq\beta _{2}\left( 1-\beta _{1}\right)
f_{n+1}f_{n+2}.\medskip $

From the above, we proved\medskip

\textbf{Theorem 2.5.} \textit{The norm of the} $n$th \textit{generalized} 
\textit{Fibonacci quaternion} $H_{n}^{p,q}$ \textit{in a generalized
quaternion algebra is}%
\begin{equation*}
\boldsymbol{n}\left( H_{n}^{p,q}\right) \text{=}p^{2}h_{2n}^{1+2\beta
_{2},3\beta _{2}}\text{+}p^{2}(\beta _{1}\text{-}1)h_{2n+1}^{1+2\beta
_{2},\beta _{2}}\text{+}q^{2}h_{2n+2}^{1+2\beta _{2},3\beta _{2}}\text{+}%
q^{2}(\beta _{1}\text{-}1)h_{2n+3}^{1+2\beta _{2},\beta _{2}}\text{-}\newline
\end{equation*}%
\begin{equation*}
-2p\left( \beta _{1}-1\right) \left( p\beta _{2}+p+q\right)
f_{n-1}f_{n}-2q^{2}(\beta _{1}-1)\left( 1\text{+}\beta _{2}\right)
f_{n}f_{n+1}+\newline
\end{equation*}%
\begin{equation}
+h_{2n+1}^{2pq\beta _{1},2pq\beta _{1}\beta _{2}}+2pq\beta _{1}\beta
_{2}(f_{2n}+f_{2n+3})+2pq\beta _{2}\left( 1-\beta _{1}\right) f_{n+1}f_{n+2}.
\tag{2.8}
\end{equation}%
$\Box $

It is known that the expression for the $n$th term of a Fibonacci element is 
\begin{equation}
f_{n}=\frac{1}{\sqrt{5}}[\alpha ^{n}-\beta ^{n}]=\frac{\alpha ^{n}}{\sqrt{5}}%
[1-\frac{\beta ^{n}}{\alpha ^{n}}],  \tag{2.9}
\end{equation}%
where $\alpha =\frac{1+\sqrt{5}}{2}$ and $\beta =\frac{1-\sqrt{5}}{2}.$

From the above, we can compute the following 
\begin{equation*}
\underset{n\rightarrow \infty }{\lim }\boldsymbol{n}\left( F_{n}\right) =%
\underset{n\rightarrow \infty }{\lim }(f_{n}^{2}+\beta _{1}f_{n+1}^{2}+\beta
_{2}f_{n+2}^{2}+\beta _{1}\beta _{2}f_{n+3}^{2})=
\end{equation*}%
\begin{equation*}
\underset{n\rightarrow \infty }{=\underset{n\rightarrow \infty }{\lim }(%
\frac{\alpha ^{2n}}{5}\text{+}\beta _{1}\frac{\alpha ^{2n+2}}{5}\text{+}%
\beta _{2}\frac{\alpha ^{2n+4}}{5}\text{+}\beta _{1}\beta _{2}\frac{\alpha
^{2n+6}}{5})=}
\end{equation*}%
\begin{equation*}
=sgnE(\beta _{1},\beta _{2})\cdot \infty
\end{equation*}%
where\newline
$E(\beta _{1},\beta _{2})=(\frac{1}{5}+\frac{\beta _{1}}{5}\alpha ^{2}+\frac{%
\beta _{2}}{5}\alpha ^{4}+\frac{\beta _{1}\beta _{2}}{5}\alpha ^{6})=$%
\newline
$=\frac{1}{5}\left( 1+\beta _{1}\left( \alpha +1\right) +\beta _{2}\left(
3\alpha +2\right) +\beta _{1}\beta _{2}\left( 8\alpha +5\right) \right) =$%
\newline
$=\frac{1}{5}[1+\beta _{1}+2\beta _{2}+5\beta _{1}\beta _{2}+\alpha \left(
\beta _{1}+3\beta _{2}+8\beta _{1}\beta _{2}\right) ],$ since $\alpha
^{2}=\alpha +1.$

If \ $E(\beta _{1},\beta _{2})>0,$ there exist a number $n_{1}\in \mathbb{N}$
such that for all\newline
$n\geq n_{1}$ we have 
\begin{equation*}
h_{2n+2}^{1+2\beta _{2},3\beta _{2}}+(\beta _{1}-1)h_{2n+3}^{1+2\beta
_{2},\beta _{2}}-2(\beta _{1}-1)\left( 1+\beta _{2}\right) f_{n}f_{n+1}>0.
\end{equation*}%
In the same way, if $E(\beta _{1},\beta _{2})<0,$ there exist a number $%
n_{2}\in \mathbb{N}$ such that for all $n\geq n_{2}$ we have 
\begin{equation*}
h_{2n+2}^{1+2\beta _{2},3\beta _{2}}+(\beta _{1}-1)h_{2n+3}^{1+2\beta
_{2},\beta _{2}}-2(\beta _{1}-1)\left( 1+\beta _{2}\right) f_{n}f_{n+1}<0.
\end{equation*}%
Therefore for \ all $\beta _{1},\beta _{2}\in \mathbb{R}$ with $E(\beta
_{1},\beta _{2})\neq 0,$ in the algebra $\mathbb{H}\left( \beta _{1},\beta
_{2}\right) $ there is a natural number $n_{0}=\max \{n_{1},n_{2}\}$ such
that $\boldsymbol{n}\left( F_{n}\right) \neq 0,$ hence $F_{n}$ is an
invertible element for all $n\geq n_{0}.$ Using the same arguments, we can
compute

\bigskip

\begin{equation*}
\underset{n\rightarrow \infty }{\lim }\left( \boldsymbol{n}\left(
H_{n}^{p,q}\right) \right) =\underset{n\rightarrow \infty }{\lim }\left(
h_{n}^{2}+\beta _{1}h_{n+1}^{2}+\beta _{2}h_{n+2}^{2}+\beta _{1}\beta
_{2}h_{n+3}^{2}\right) =
\end{equation*}%
\begin{equation*}
=\underset{n\rightarrow \infty }{\lim }[\left( pf_{n-1}\text{+}qf_{n}\right)
^{2}\text{+}\beta _{1}\left( pf_{n}\text{+}qf_{n+1}\right) ^{2}\text{+}\beta
_{2}\left( pf_{n+1}\text{+}qf_{n+2}\right) ^{2}+
\end{equation*}%
\begin{equation*}
\text{+}\beta _{1}\beta _{2}\left( pf_{n+2}\text{+}qf_{n+3}\right) ^{2}]=
\end{equation*}%
\begin{equation*}
=sgnE^{\prime }(\beta _{1},\beta _{2})\cdot \infty
\end{equation*}%
where\newline
$E^{\prime }(\beta _{1},\beta _{2})=\frac{1}{5}[\left( p+\alpha q\right)
^{2}+\beta _{1}\left( p\alpha +\alpha ^{2}q\right) ^{2}+\beta _{2}\left(
p\alpha ^{2}+\alpha ^{3}q\right) ^{2}+$\newline
$+\beta _{1}\beta _{2}\left( p\alpha ^{3}+\alpha ^{4}q\right) ^{2}]=$\newline
$=\frac{1}{5}\left( p+\alpha q\right) ^{2}[1+\beta _{1}\alpha ^{2}+\beta
_{2}\alpha ^{4}+\beta _{1}\beta _{2}\alpha ^{6}]=$\newline
$=\frac{1}{5}\left( p+\alpha q\right) ^{2}E(\beta _{1},\beta _{2}).$

Therefore for \ all $\beta _{1},\beta _{2}\in \mathbb{R}$ with $E^{\prime
}(\beta _{1},\beta _{2})\neq 0$ in the algebra $\mathbb{H}\left( \beta
_{1},\beta _{2}\right) $ there exist a natural number $n_{0}^{\prime }$ such
that $\boldsymbol{n}\left( H_{n}^{p,q}\right) \neq 0,$ hence $H_{n}^{p,q}$
is an invertible element for all $n\geq n_{0}^{\prime }.$

Now, we proved\medskip \medskip

\textbf{Theorem 2.6.} \textit{For all} $\beta _{1},\beta _{2}\in \mathbb{R}$ 
\textit{with} $E^{\prime }(\beta _{1},\beta _{2})\neq 0,$ \textit{there
exist a natural number} $n^{\prime }$ \textit{such that for all} $n\geq
n^{\prime }$ \textit{\ Fibonacci elements} $\ F_{n}$ \textit{and generalized
Fibonacci elements} $H_{n}^{p,q}$ \textit{are invertible elements in the
algebra} $\mathbb{H}\left( \beta _{1},\beta _{2}\right) .\Box \medskip $

\textbf{Remark 2.7.} Algebra $\mathbb{H}\left(\beta_1,\beta _{2}\right) $ is
not always a division algebra, and sometimes can be difficult to find an
example \ of invertible element. Above Theorem \ provides us infinite sets
of invertible elements in this algebra, namely Fibonacci elements and
generalized Fibonacci elements.

\textit{\ } \ 
\begin{equation*}
\end{equation*}

{\Large 3. \ Fibonacci-Narayana Quaternions\bigskip }

\begin{equation*}
\end{equation*}

\bigskip In this section, we will study some properties of
Fibonacci-Narayana elements in the algebra $\mathbb{H}\left( \beta
_{1},\beta _{2}\right) .\medskip $

\textbf{Theorem 3.1}. \textit{For the Fibonacci-Narayana quaternion} $U_{n}$,%
\textit{\ we have}

\begin{equation*}
a)\,\sum\limits_{m=0}^{n}U_{m}=U_{n+3}-U_{2},
\end{equation*}%
\begin{equation*}
b)\,\sum\limits_{m=0}^{n}U_{3m}=U_{3n+1}-1-e_{4}.
\end{equation*}

\textbf{Proof.} a) 
\begin{equation*}
\sum\limits_{m=0}^{n}U_{m}=\sum\limits_{m=0}^{n}u_{m}+e_{2}\sum%
\limits_{m=1}^{n+1}u_{m}+e_{3}\sum\limits_{m=2}^{n+2}u_{m}+e_{4}\sum%
\limits_{m=3}^{n+3}u_{m}=
\end{equation*}%
Since $u_{0}=0,\,$\ we consider that the term $\sum\limits_{m=0}^{n}u_{m}$
is equal with $\sum\limits_{m=1}^{n}u_{m}.$ We can use property $1)$ from
Introduction and we obtain 
\begin{equation*}
=u_{n+3}-1+e_{2}(u_{n+4}-1)+e_{3}(u_{n+5}-2)+e_{4}(u_{n+6}-3)=
\end{equation*}%
\begin{equation*}
=U_{n+3}-(1+e_{2}+2e_{3}+3e_{4})=U_{n+3}-U_{2}.
\end{equation*}

b)$\,\,$Since $u_{0}=0,$ the term $\sum\limits_{m=0}^{n}u_{3m}$ is equal
with $\sum\limits_{m=1}^{n}u_{3m}$, therefore 
\begin{equation*}
\sum\limits_{m=0}^{n}U_{3m}=\sum\limits_{m=0}^{n}u_{3m}+e_{2}\sum%
\limits_{m=0}^{n}u_{3m+1}+e_{3}\sum\limits_{m=0}^{n}u_{3m+2}+e_{4}\sum%
\limits_{m=0}^{n}u_{3m+3}=
\end{equation*}%
using properties $4)$, $2)$, $3)$, and again $4)$, we have 
\begin{equation*}
=u_{3n+1}-1+u_{3n+2}e_{2}+u_{3n+3}e_{3}+(u_{3n+4}-1)e_{4}=U_{3n+1}-1-e_{4}.
\end{equation*}%
$\Box \medskip $

Let $\{u_{n}\}$ be a Fibonacci-Narayana sequence, and let $U_{n}=u_{n}\cdot
1+u_{n+1}e_{2}+u_{n+2}e_{3}+u_{n+3}e_{4}$ be the $n$th Fibonacci-Narayana
quaternion.

The function $f(x)=a_{0}+a_{1}x+a_{2}x^{2}+\ldots +a_{n}x^{n}+\ldots $ is
called\textit{\ the generating function} for the sequence $%
\{a_{0},a_{1},a_{2},\ldots \}$. In [Ha; 12], the author found a generating
function for Fibonacci quaternions. In the following theorem, we established
the generating function for Fibonacci-Narayana quaternions.\medskip\ 

\textbf{Theorem 3.2.} \textit{The generating function for the
Fibonacci-Narayana quaternion} $U_{n}$ \textit{is} 
\begin{equation}
G(t)\text{=}\frac{U_{0}\text{+}(U_{1}\text{-}U_{0})t\text{+}(U_{2}\text{-}%
U_{1})t^{2}}{1-t-t^{3}}\text{=}\frac{e_{1}\text{+}e_{2}\text{+}e_{3}\text{+}%
(1\text{+}e_{3})t+(e_{2}\text{+}e_{3})t^{2}}{1-t-t^{3}}.  \tag{3.1}
\end{equation}

\textbf{Proof.} Assuming that the generating function of the \ quaternion%
\newline
Fibonacci-Narayana sequence $\{U_{n}\}$ has the form $G(t)=\sum%
\limits_{n=0}^{\infty }U_{n}\,t^{n}$, we obtain that\newline
$\sum\limits_{n=0}^{\infty }U_{n}t^{n}-t\sum\limits_{n=0}^{\infty
}U_{n}t^{n}-t^{3}\sum\limits_{n=0}^{\infty }U_{n}t^{n}=$\newline
$=U_{0}+U_{1}t+U_{2}t^{2}+U_{3}t^{3}+...-$\newline
$-U_{0}t-U_{1}t^{2}-U_{2}t^{3}-U_{3}t^{4}-...-$\newline
$-U_{0}t^{3}-U_{1}t^{4}-U_{2}t^{5}-U_{3}t^{6}-....=$\newline
$=U_{0}$+$(U_{1}-U_{0})t$+$(U_{2}-U_{1})t^{2},$ since $%
U_{n}=U_{n-1}+U_{n-3},n\geq 3$ and the coefficients of $t^{n}$ for $%
\,\,n\geq 3$ are equal with zero.

It results%
\begin{equation*}
U_{0}\text{+}(U_{1}-U_{0})t\text{+}(U_{2}-U_{1})t^{2}\text{=}%
\sum\limits_{n=0}^{\infty }U_{n}t^{n}\,(1-t-t^{3}),
\end{equation*}%
or in equivalent form 
\begin{equation*}
\frac{U_{0}\text{+}(U_{1}-U_{0})t\text{+}(U_{2}-U_{1})t^{2}}{1-t-t^{3}}\text{%
=}\sum\limits_{n=0}^{\infty }U_{n}t^{n}.
\end{equation*}%
The theorem is proved. $\Box \medskip $

\textbf{Theorem 3.3.} (The Binet-Cauchy formula for Fibonacci-Narayana
numbers) \textit{Let} $u_{n}=u_{n-1}+u_{n-3},n\geq 3$ \textit{be the} $n$%
\textit{th} \textit{Fibonacci-Narayana number, then }%
\begin{equation}
u_{n}\text{=}\frac{1}{\left( \alpha \text{-}\beta \right) (\beta \text{-}%
\gamma )\left( \gamma \text{-}\alpha \right) }\left[ \alpha ^{n+1}\left(
\gamma \text{-}\beta \right) \text{+}\beta ^{n+1}\left( \alpha \text{-}%
\gamma \right) \text{+}\gamma ^{n+1}\left( \beta \text{-}\alpha \right) %
\right] ,  \tag{3.2}
\end{equation}%
\textit{where} $\alpha ,\beta ,\gamma \,$\ \textit{are the solutions of the
equation} $t^{3}-t^{2}-1=0.\medskip $

\textbf{Proof.} Supposing that $u_{n}=A\alpha ^{n}+B\beta ^{n}+C\gamma
^{n},A,B,C\in \mathbb{C}$ and using the recurrence formula for the
Fibonacci-Narayana numbers, $u_{n}=u_{n-1}+u_{n-3}$, it results that $\alpha
,\beta ,\gamma \,$\ are the solutions of the equation $t^{3}-t^{2}-1=0.$
Since $u_{0}=0,u_{1}=1,u_{2}=1,$ we obtain the following system%
\begin{equation}
\left\{ 
\begin{array}{c}
A+B+C=0, \\ 
A\alpha +B\beta +C\gamma =1, \\ 
A\alpha ^{2}+B\beta ^{2}+C\gamma ^{2}=1.%
\end{array}%
\right.  \tag{3.3}
\end{equation}%
The determinant of this system is a Vandermonde determinant and can be
computed easily. It is $\Delta =\left( \alpha \text{-}\beta \right) (\beta $-%
$\gamma )\left( \gamma \text{-}\alpha \right) \neq 0.$

Using the Cramer's rule, the solutions of the system $\left( 3.3\right) $ are%
\newline
\begin{equation*}
A=\frac{\alpha \left( \gamma -\beta \right) }{\left( \alpha -\beta \right)
(\beta -\gamma )\left( \gamma -\alpha \right) }=\frac{\alpha }{\left( \beta
-\alpha \right) \left( \gamma -\alpha \right) }\medskip,
\end{equation*}
\begin{equation*}
B=\frac{\beta \left( \alpha -\gamma \right) }{\left( \alpha -\beta \right)
(\beta -\gamma )\left( \gamma -\alpha \right) }=\frac{\beta }{\left( \alpha
-\beta \right) (\gamma -\beta )}\medskip,
\end{equation*}
\begin{equation*}
C=\frac{\gamma \left( \beta -\alpha \right) }{\left( \alpha -\beta \right)
(\beta -\gamma )\left( \gamma -\alpha \right) }=\frac{\gamma }{(\beta
-\gamma )\left( \alpha -\gamma \right) },
\end{equation*}
therefore relation $\left( 3.2\right) $ is true.$~\Box \medskip $

\textbf{Theorem 3.4.} (The Binet-Cauchy formula for the Fibonacci-Narayana
quaternions) \textit{Let} $U_{n}=u_{n}\cdot
1+u_{n+1}e_{2}+u_{n+2}e_{3}+u_{n+3}e_{4}$ \textit{be the} $n$\textit{th
Fibonacci-Narayana quaternion, then} 
\begin{equation}
U_{n}\text{=}D\frac{\alpha ^{n+1}}{\left( \beta \text{-}\alpha \right)
\left( \gamma \text{-}\alpha \right) }\text{+}E\frac{\beta ^{n+1}}{\left(
\alpha \text{-}\beta \right) (\gamma \text{-}\beta )}\text{+}F\frac{\gamma
^{n+1}}{(\beta \text{-}\gamma )\left( \alpha \text{-}\gamma \right) }, 
\tag{3.4}
\end{equation}%
\textit{where} $\alpha ,\beta ,\gamma \,$\ \textit{are the solutions of the
equation} $t^{3}-t^{2}-1=0$ \textit{and} \newline
\begin{equation*}
D=1+\alpha e_{1}+\alpha ^{2}e_{2}+\alpha ^{3}e_{3},
\end{equation*}
\begin{equation*}
E=1+\beta e_{1}+\beta ^{2}e_{2}+\beta ^{3}e_{3},
\end{equation*}
\begin{equation*}
F=1+\gamma e_{1}+\gamma ^{2}e_{2}+\gamma ^{3}e_{3}.\medskip
\end{equation*}

\textbf{Proof.} Using relation $\left( 3.2\right) ,~$we have that\newline
$U_{n}$=$u_{n}\cdot 1$+$u_{n+1}e_{2}$+$u_{n+2}e_{3}$+$u_{n+3}e_{4}$=$%
\medskip \newline
$=$\frac{1}{\left( \alpha \text{-}\beta \right) (\beta \text{-}\gamma
)\left( \gamma \text{-}\alpha \right) }[\left( \alpha ^{n\text{+}1}\left(
\gamma \text{-}\beta \right) \text{+}\beta ^{n\text{+}1}\left( \alpha \text{-%
}\gamma \right) \text{+}\gamma ^{n\text{+}1}\left( \beta \text{-}\alpha
\right) \right) \cdot 1$+\newline
+$\left( \alpha ^{n+2}\left( \gamma -\beta \right) +\beta ^{n+2}\left(
\alpha -\gamma \right) +\gamma ^{n+2}\left( \beta -\alpha \right) \right)
e_{1}$+\newline
+$\left( \alpha ^{n+3}\left( \gamma -\beta \right) +\beta ^{n+3}\left(
\alpha -\gamma \right) +\gamma ^{n+3}\left( \beta -\alpha \right) \right)
e_{2}$+\newline
$+\left( \alpha ^{n+4}\left( \gamma -\beta \right) +\beta ^{n+4}\left(
\alpha -\gamma \right) +\gamma ^{n+4}\left( \beta -\alpha \right) \right)
e_{3}]$=\newline
=$\frac{1}{\left( \alpha -\beta \right) (\beta -\gamma )\left( \gamma
-\alpha \right) }[\alpha ^{n+1}\left( \gamma \text{-}\beta \right) \left( 1%
\text{+}\alpha e_{1}\text{+}\alpha ^{2}e_{2}\text{+}\alpha ^{3}e_{3}\right) $%
+\newline
$+\beta ^{n+1}\left( \alpha \text{-}\gamma \right) \left( 1+\beta
e_{1}+\beta ^{2}e_{2}+\beta ^{3}e_{3}\right) +\newline
+\gamma ^{n+1}\left( \beta \text{-}\alpha \right) \left( 1+\gamma
e_{1}+\gamma ^{2}e_{2}+\gamma ^{3}e_{3}\right) ].$\newline
$\Box \medskip $

For negative $n$, the $n$th Fibonacci-Narayana number will be defined as $%
u_{n}=u_{n+3}-u_{n+2},\,\,u_{0}=0,u_{1}=1,u_{2}=1$. Accordingly defined the
Fibonacci-Narayana quaternion $U_{n}$ for negative $n$.\medskip

\textbf{Theorem 3.5.} \textit{Let} $U_{n}=u_{n}\cdot
1+u_{n+1}e_{2}+u_{n+2}e_{3}+u_{n+3}e_{4}$ \textit{be the} $n$\textit{th
Fibonacci-Narayana quaternion, therefore the following relations are true:}

1) $\underset{i=0}{\overset{n}{\sum }}\complement
_{n}^{i}U_{2n-2i-1}=U_{3n-1}.\medskip $

2)$\underset{i=0}{\overset{n}{\sum }}\complement
_{n}^{i}U_{3n-2i-1}=U_{4n-1.}\medskip $

\textbf{Proof.}

1) \ Using the Newton's formula, it results that\medskip \newline
$\left( t^{2}+1\right) ^{n}=\complement _{n}^{0}\left( t^{2}\right)
^{n}+\complement _{n}^{1}\left( t^{2}\right) ^{n-1}+\complement
_{n}^{2}\left( t^{2}\right) ^{n-2}+...+\complement _{n}^{n}=\medskip $%
\newline
=$\complement _{n}^{0}t^{2n}+\complement _{n}^{1}t^{2n-2}+\complement
_{n}^{2}t^{2n-4}+...+\complement _{n}^{n}.$ From here, we have that\medskip 
\newline
$\underset{i=0}{\overset{n}{\sum }}\complement _{n}^{i}U_{2n-2i-1}$=$%
\complement _{n}^{0}U_{2n-1}+\complement _{n}^{1}U_{2n-3}+\complement
_{n}^{2}U_{2n-5}+...+$ $\complement _{n}^{n}U_{-1}=\medskip $\newline
=$\complement _{n}^{0}\left( D\frac{\alpha ^{2n}}{\left( \beta \text{-}%
\alpha \right) \left( \gamma \text{-}\alpha \right) }+E\frac{\beta ^{2n}}{%
\left( \alpha \text{-}\beta \right) (\gamma \text{-}\beta )}+F\frac{\gamma
^{2n}}{(\beta \text{-}\gamma )\left( \alpha \text{-}\gamma \right) }\right)
+\medskip $\newline
+$\complement _{n}^{1}\left( D\frac{\alpha ^{2n-2}}{\left( \beta \text{-}%
\alpha \right) \left( \gamma \text{-}\alpha \right) }+E\frac{\beta ^{2n-2}}{%
\left( \alpha \text{-}\beta \right) (\gamma \text{-}\beta )}+F\frac{\gamma
^{2n-2}}{(\beta \text{-}\gamma )\left( \alpha \text{-}\gamma \right) }%
\right) +...+\medskip $\newline
+$\complement _{n}^{n}\left( D\frac{1}{\left( \beta \text{-}\alpha \right)
\left( \gamma \text{-}\alpha \right) }+E\frac{1}{\left( \alpha \text{-}\beta
\right) (\gamma \text{-}\beta )}+F\frac{1}{(\beta \text{-}\gamma )\left(
\alpha \text{-}\gamma \right) }\right) =\medskip $\newline
=$D\frac{1}{\left( \beta \text{-}\alpha \right) \left( \gamma \text{-}\alpha
\right) }(\complement _{n}^{0}\alpha ^{2n}+\complement _{n}^{1}\alpha
^{2n-2}+...+\complement _{n}^{n}1)+\medskip $\newline
+$E\frac{1}{\left( \alpha \text{-}\beta \right) (\gamma \text{-}\beta )}%
(\complement _{n}^{0}\beta ^{2n}+\complement _{n}^{1}\beta
^{2n-2}+...+\complement _{n}^{n}1)+\medskip $\newline
+$F\frac{1}{(\beta \text{-}\gamma )\left( \alpha \text{-}\gamma \right) }%
(\complement _{n}^{0}\gamma ^{2n}+\complement _{n}^{1}\gamma
^{2n-2}+...\complement _{n}^{n}1)=\medskip $\newline
=$D\frac{1}{\left( \beta \text{-}\alpha \right) \left( \gamma \text{-}\alpha
\right) }\left( \alpha ^{2}\text{+}1\right) ^{n}$+$E\frac{1}{\left( \alpha 
\text{-}\beta \right) (\gamma \text{-}\beta )}\left( \beta ^{2}\text{+}%
1\right) ^{n}$+$F\frac{1}{(\beta \text{-}\gamma )\left( \alpha \text{-}%
\gamma \right) }\left( \gamma ^{2}\text{+}1\right) ^{n}=\medskip $\newline
=$D\frac{1}{\left( \beta \text{-}\alpha \right) \left( \gamma \text{-}\alpha
\right) }\alpha ^{3n}+E\frac{1}{\left( \alpha \text{-}\beta \right) (\gamma 
\text{-}\beta )}\beta ^{3n}+F\frac{1}{(\beta \text{-}\gamma )\left( \alpha 
\text{-}\gamma \right) }\gamma ^{3n}=U_{3n-1}.\medskip \medskip $\newline
We used that $\alpha ^{3}=\alpha ^{2}+1,\beta ^{3}=\beta ^{2}+1,\gamma
^{3}=\gamma ^{2}+1.$

2) Since $t^{3}=t^{2}+1$, starting from relation $\left( t^{3}+t\right)
^{n}=t^{4n}$, for\newline
$t\in \{\alpha ,\beta ,\gamma \},$ by straightforward calculations as in 2),
we obtain the asked relation. $\Box \medskip $

\begin{equation*}
\end{equation*}

\textbf{Conclusions.} In this paper we investigated some new properties of
generalized Fibonacci quaternions and Fibonacci-Narayana quaternions. Since
Fibonacci-Narayana quaternions was not intensive studied until now, we
expect to find in the future more and surprising new properties. We study
these elements for the beauty of the relations obtained, but the main reason
\ is that the elements \ of \ this type, namely \ Fibonacci $X$ elements,
where $X\in \{quaternions,generalized\ quaternions\},$ can provide us many
important information in the algebra $\mathbb{H}\left( \beta _{1},\beta
_{2}\right) ,$ as for example: sets of invertible elements in algebraic
structures without division. 
\begin{equation*}
\end{equation*}

\textbf{Acknowledgements.} Authors thank referee for his/her patience and
suggestions which help us to improve this paper.

\begin{equation*}
\end{equation*}%
\textbf{References}%
\begin{equation*}
\end{equation*}

[Al, Jo; 96 ] J.P. Allouche, \ T. Johnson, \textit{Nrayana's Cows and \
Delayed Morphisms}, recherche.ircam.fr/equipes/repmus/jim96/actes/Allouche.ps

[Cu; 76] I. Cucurezeanu, \ \textit{Exercices in Arithmetics and Numbers
Theory}, Ed. Tehnic\u{a}, Bucharest, 1976. [in Romanian]

[Di, St; 03] T.~V.~Didkivska, M.~V.~St'opochkina, \textit{Properties of
Fibonacci-Narayana numbers}, In the World of Mathematics, \textbf{9}%
(1)(2003), 29--36. [in Ukrainian]

[Dr, Gi, Gr, Wa; 03], A. Dress, R. Giegerich, S. Gr\H{u}newald, H. Wagner, \ 
\textit{Fibonacci-Cayley Numbers and Repetition Patterns in Genomic DNA},
Ann. Comb. \textbf{7}(2003) 259--279.

[Ha; 12] S. Halici, \textit{On Fibonacci Quaternions}, Adv. in Appl.
Clifford Algebras, \textbf{22}(2)(2012), 321-327.

[Ho; 61] A. F. Horadam, \textit{A Generalized Fibonacci Sequence}, Amer.
Math. Monthly, \textbf{68}(1961), 455-459.

[Ho; 63] A. F. Horadam, \textit{Complex Fibonacci Numbers and Fibonacci
Quaternions}, Amer. Math. Monthly, \textbf{70}(1963), 289-291.

[Ka; 04] S. Kak, \textit{The Golden Mean and the Physics of Aesthetics},
Archive of Physics: physics/0411195, 2004,\newline
http://arxiv.org/ftp/physics/papers/0411/0411195.pdf.

[Ko; 01] T.~Koshy, \textit{Fibonacci and Lucas Numbers with Applications}, A
Wiley-Interscience publication, U.S.A, 2001.

[Pa; 85] P. Singh, \textit{The So-called Fibonacci Numbers in Ancient and
Medieval India}, Historia Mathematica, \textbf{12}(1985), 229-244.

[Sa-Mu; 82] P. V. Satyanarayana Murthy, \textit{Fibonacci-Cayley Numbers},
The Fibonacci Quarterly, \textbf{20}(1)(1982), 59-64.

[Si; 36] A.N. Singh, \textit{On the use of series in Hindu mathematics},
Osiris 1, 1936, \ p. 606-628. \newline
(http://www.anaphoria.com/hindu.PDF)

[Si; 85] P. Singh, \textit{The so-called Fibonacci numbers in ancient and
medieval India}, Historia Mathematica \textbf{12}(1985), 229-244.

[Sh; 06] V.~Shpakivskyi, \textit{On the one property of Fibonacci-Narayana
numbers}, Mathematics, no. \textbf{11}(359)(2006), 12--14. [in Ukrainian]

[Sw; 73] M. N. S. Swamy, \textit{On generalized Fibonacci Quaternions}, The
Fibonacci Quaterly, \textbf{11(}5\textbf{)}(1973), 547-549.

[Yu; 61] A.~P.~Yushkevich, \textit{History of Mathematics in the Middle Ages}%
, Gos. izd. fiz.-mat. lit., Moskow, 1961.%
\begin{equation*}
\end{equation*}

Cristina FLAUT

{\small Faculty of Mathematics and Computer Science,}

{\small Ovidius University,}

{\small Bd. Mamaia 124, 900527, CONSTANTA,}

{\small ROMANIA}

{\small http://cristinaflaut.wikispaces.com/}

{\small http://www.univ-ovidius.ro/math/}

{\small e-mail:}

{\small cflaut@univ-ovidius.ro}

{\small cristina\_flaut@yahoo.com}%
\begin{equation*}
\end{equation*}

Vitalii \ SHPAKIVSKYI

{\small Department of Complex Analysis and Potential Theory}

{\small \ Institute of Mathematics of the National Academy of Sciences of
Ukraine,}

{\small \ 3, Tereshchenkivs'ka st.}

{\small \ 01601 Kiev-4}

{\small \ UKRAINE}

{\small \ http://www.imath.kiev.ua/\symbol{126}complex/}

{\small \ e-mail: shpakivskyi@mail.ru}

\end{document}